# DESIGN OF LOW-THRUST MULTI-GRAVITY ASSIST TRAJECTORIES TO EUROPA

**Massimiliano VASILE[1], Stefano CAMPAGNOLA[2]**

[1]*Department of Aerospace Engineering,
University of Glasgow, James Watt Building G12 8QQ, Glasgow, UK
m.vasile@aero.gla.ac.uk*

[2]*University of Southern California
AME Department - Aerospace and Mechanical Engineering, Los Angeles, California, USA*

**ABSTRACT :**

This paper presents the design of a mission to Europa using solar electric propulsion as main source of thrust. A direct transcription method based on Finite Elements in Time was used for the design and optimisation of the entire low-thrust gravity assist transfer from the Earth to Europa. Prior to that, a global search algorithm was used to generate a set of suitable first guess solutions for the transfer to Jupiter, and for the capture in the Jovian system. In particular, a fast deterministic search algorithm was developed to find the most promising set of swing-bys to reach Jupiter. A second fast search algorithm was developed to find the best sequence of swing-bys of the Jovian moons. After introducing the global search algorithms and the direct transcription through Finite Elements in Time, the paper presents a number of first guess solutions and a fully optimised transfer from the Earth to Europa.

## 1 - INTRODUCTION

Europa, one of the moons of Jupiter, is among the most interesting targets in the solar system for exobiology. The presence of water, probably liquid under a superficial layer of ice, suggests the possibility of prebiotic life. A mission aimed at the exploration of this peculiar moon would be extremely interesting but, at the same time, quite challenging due to the distance from both the Earth and Sun. Such a mission would be extremely demanding in terms of communications, power and *Δv*. All the missions to Jupiter, or more generally to the outer part of the solar system, such as Galileo [1] or Cassini [2], have employed chemical propulsion as the main propulsion system and radioisotope thermal generators (RTGs) to generate the required power on-board. An alternative would be to use solar electric propulsion or a combination of electric and chemical propulsion.

Recent studies have explored the possibility of missions to the icy moons of Jupiter based on nuclear electric propulsion. Parcher et al. [3][4] proposed a number of transfer trajectories to Jupiter optimising both time of flight and propellant mass. Sims presented an analysis of a potential mission to the icy moons, with some results on the optimal design of a transfer to Jupiter and of a tour of the Jovian moons [5]. All these studies were addressing the design of the Jupiter Icy Moons Orbiter (JIMO). The main target was Europa due to the aforementioned-suspected ocean where simple alien life is a possibility in our solar system. Ganymede and Callisto, which are now thought to have liquid, salty oceans beneath their icy surfaces, were included in the targets of the JIMO mission. The spacecraft was planned to be propelled by eight ion engines and powered by a small fission reactor.

In this paper, we investigate the design of a solar electric propulsion (SEP) mission to Jupiter's moons. The requirements in terms of $\Delta v$ are reduced by using multiple Gravity Assist (GA) manoeuvres at Venus, Earth, Mars, and the Jovian moons. The inner planets and the Jovian moons are modelled as simple point masses with no gravity field, while Jupiter is modelled as a full gravitational body. The main goal of this study is to assess the actual feasibility of a mission to Europa using solar electric propulsion as main source of thrust, and to investigate the related issues.

The design of a transfer trajectory combining SEP and GA can be regarded as a general trajectory optimisation problem. The dynamics of the spacecraft is mainly governed by the gravitational attraction of the Sun when the spacecraft is outside the sphere of influence of a planet, and by the gravitational attraction of the planet during a gravity assist manoeuvre. Low-thrust propulsion is used to shape the trajectory arcs between two subsequent encounters and to obtain the best incoming conditions for a swing-by. In this work, the trajectory was split into several phases, with each phase corresponding to a trajectory arc connecting two celestial bodies. For each phase, a Finite Elements in Time technique [6][7] was used to transcribe the differential equations governing the dynamics of the spacecraft into a set of algebraic nonlinear equations. Direct Finite Elements Transcription (DFET) has already been proven to be an effective direct method for the design of very complex trajectories with multiple gravity assist manoeuvres and low-thrust arcs, including the optimisation of system parameters and a full *n*-body dynamics for gravity assist manoeuvres [6][7].

Prior to the optimisation of the low-thrust trajectory, three algorithms, implementing three reduced trajectory models of increasing complexity, were used to generate first guess solutions for both the transfer from the Earth to the Jupiter, and for the tour of the moons. Of the three algorithms, two were used to explore the search space for all possible solutions for a transfer to Jupiter, and one was used to design the tour of Jupiter's moons. In particular, they were used to automatically find a set of potentially optimal swing-by sequences (i.e. the number of gravity assist manoeuvres and the celestial body of each manoeuvre). Unlike other approaches, which implement a reduced model for low-thrust trajectories [8][9], here we use a reduced model for a multi-impulse trajectory to generate a first guess for a low-thrust trajectory.

The two algorithms for the design of the transfer to Jupiter implement a trajectory model conceptually similar to S-TOUR [10][11][12], and a systematic search. The design is a two-step process; the first step identifies the sequence and a possible interval of launch dates, and the second computes an accurate multi-impulsive trajectory. The first step makes use of a very simple model, which allows for the computation of a very large number of potential solutions in a matter of seconds. The third algorithm was devised to generate potentially optimal sequences of swing-bys of Jupiter's moon. The algorithm provides results that are consistent with the work of Heaton et al. [13] with a very low computational effort.

The paper is structured as follows: after a description of the dynamic models, the three algorithms for the generation of first guesses are presented in detail. A description of the DFET method then follows. Finally, prior to the results section, we present an analysis of the capture manoeuvre at Jupiter using a gravity assist manoeuvre of one of the moons.

## 2 - TRAJECTORY MODEL

The whole transfer to Europa was decomposed into two main parts: an interplanetary cruise from the Earth to Jupiter and a capture in the Jovian system with a tour of the Jovian moons. Along the first part of the transfer the spacecraft is subject to the gravity attraction of the Sun, the gravity disturbance of Jupiter and the thrust of the electric propulsion engine. During the capture part of the transfer the spacecraft is subject to the gravity attraction of Jupiter, the gravity disturbance of the Sun and the thrust of the electric engine. The maximum thrust delivered by the electric engine was computed as a function of the power provided by the solar panels. For each one of the two parts, a different reference system was considered. Each part then was split into a number of phases. Each phase corresponds to a transfer arc connecting two planets or two moons.

An appropriate set of inter-phase constraints was used to assembled together all the phases, forming a single Nonlinear Programming problem comprising the two parts of the transfer. The NLP problem was then solved with a sparse sequential quadratic programming algorithm.

### 2.1 - DYNAMIC MODELS

The spacecraft is modelled as a point mass subject to the gravitational attraction of the Sun and Jupiter, and to the thrust delivered by one or more low-thrust engines. The motion of the spacecraft is described in the J2000 mean ecliptic reference frame centred on the Sun (see Figure 1) during cruise, and in the J2000 mean equatorial reference frame centred on Jupiter during capture and the tour of the Jovian moons. The three components of the thrust vector **u** represent the control:

$$\dot{\mathbf{r}} = \mathbf{v}$$
$$\dot{\mathbf{v}} = \nabla U(\mathbf{r}) + \nabla U_B(\mathbf{r}) + \frac{\mathbf{u}}{m} \quad (1)$$

where $U$ is the gravitational potential of the principal attracting body, either the Sun or Jupiter, and is a function of the position vector $\mathbf{r} = \{r_x, r_y, r_z\}^T$:

$$U(\mathbf{r}) = \frac{\mu}{|\mathbf{r}|} \quad (2)$$

The disturbing potential due to the gravity of a third body is given by:

$$U_B(\mathbf{r}) = \mu_B \left( \frac{1}{\mathbf{d}} - \frac{\langle \mathbf{d}, \rho \rangle}{\rho^3} \right) \quad (3)$$

where $\rho$ is the position vector of the perturbing body with respect to the principal one, $\mathbf{d} = \mathbf{r} - \rho$ is the position vector of the spacecraft with respect to the perturbing body and $\mu_B$ is the gravitational constant of the perturbing body. The state and the control vectors are then defined as follows:

$$\mathbf{x} = \{r_x, r_y, r_z, v_x, v_y, v_z, m\}^T \qquad \mathbf{u} = \{u_x, u_y, u_z\}^T \quad (4)$$

where $m$ is the mass of the spacecraft.

An upper bound $T_{max}$ and a lower bound $T_{min}$ was put on the thrust magnitude:

$$T_{min} \leq u = \sqrt{u_x^2 + u_y^2 + u_z^2} \leq T_{max} \quad (5)$$

The upper boundary is set equal to the maximum level of thrust provided by the selected low-thrust engine, while the lower limit was set to $1 \times 10^{-4}$ times $T_{max}$ to avoid singularities in the Hessian matrix when minimum mass problems are solved. The mass flow rate of the engine is given by,

$$\dot{m} = -\frac{u}{I_{sp} g_0} \quad (6)$$

where $I_{sp}$ is the specific impulse of the engine and $g_0$ the gravity on Earth surface.

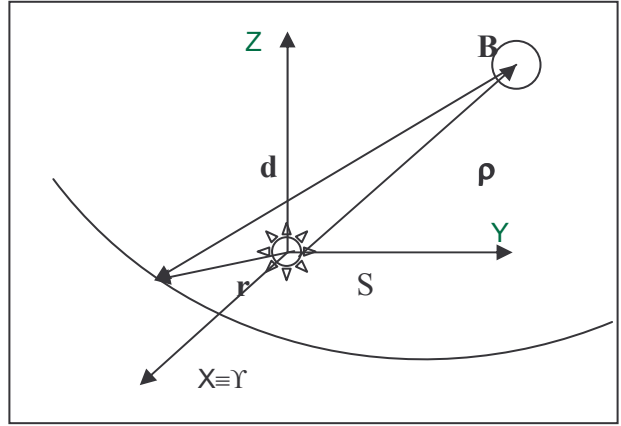

Figure 1. Heliocentric inertial reference frame, where the ecliptic plane lies in the *x-y* plane.

## 2.2 - SWING-BY MODEL

The simplest way to model a gravity assist manoeuvre is to use a linked-conic approximation where the sphere of influence of a planet is assumed to have zero radius, and the gravity manoeuvre is considered as instantaneous. Therefore, the instantaneous position vector is not affected by the swing-by:

$$\mathbf{r}_i = \mathbf{r}_0 = \mathbf{r}_P \qquad (7)$$

where $\mathbf{r}_i$ is the incoming heliocentric position, $\mathbf{r}_0$ is the outgoing heliocentric position vector and $\mathbf{r}_p$ is the planet position vector, all taken at the epoch of the encounter. For an ideal hyperbolic orbit not subject to perturbations or $\Delta v$ manoeuvres, the modulus of the incoming relative velocity $\tilde{\mathbf{v}}^-$ must be equal to the modulus of the outgoing relative velocity $\tilde{\mathbf{v}}^+$:

$$\tilde{v} = \tilde{v}^- = \tilde{v}^+ \qquad (8)$$

Furthermore, the outgoing relative velocity vector is rotated, due to gravity, by an angle $\beta$ with respect to the incoming velocity vector. Therefore, the following relation must hold:

$$\langle \tilde{\mathbf{v}}^+, \tilde{\mathbf{v}}^- \rangle = \tilde{v}^2 \cos \beta \qquad (9)$$

where the angle of rotation of the velocity vector is:

$$\beta = 2\arcsin\left(\frac{\tilde{\mu}}{\tilde{v}^2 \tilde{r}_p + \tilde{\mu}}\right) \qquad (10)$$

All quantities with a tilde are relative to the swing-by planet and $\tilde{r}_p$ is the periapsis radius of the swing-by hyperbola.

## 2.3 - POWER MODEL

One of the major issues for an interplanetary transfer using solar electric propulsion is the power available to the engine. More specifically, the thrust level depends on the power provided by the solar arrays, which, in turn, is a function of the distance from the Sun. The maximum thrust that the engine can deliver depends on the engine characteristics, such as the specific thrust $F_{sp}$ on the effective input power $P_{in}$ provided by the power system and on an efficiency coefficient $\eta_e$:

$$F_{max} = \eta_e P_{in} F_{sp} \qquad (11)$$

The effective input power is the effective power produced by the solar arrays minus the power required by all the subsystems onboard the spacecraft $P_{SS}$:

$$P^*_{in} = P_{eff} - P_{SS} \qquad (12)$$

The power delivered by the solar arrays is a function of the temperature and of the distance from the Sun:

$$P_{eff} = \eta_S \frac{P_{1AU}}{R_S^2}[1 - C_T(T_S - T_0)]\cos\alpha_{ss} \qquad (13)$$

where $P_{1AU}$ is the power generated at one Astronomical Unit (AU), $T_s$ is the temperature of solar arrays, $R_S$ is the distance from the Sun, $T_0$ the reference temperature, $C_T$ is the temperature coefficient which express the variation of performance of the array with temperature, $\eta_S$ is the efficiency of the power system and $\alpha_{ss}$ is the solar array Sun aspect angle, i.e. the angle between the normal to the cell surface and the Sun direction.

The steady state surface temperature of the solar panels is a function of the distance from the Sun:

$$T_S = \left[\frac{S_0 \alpha_s \cos\alpha_{ss}}{R_S^2 \sigma \kappa \varepsilon_s}\right]^{0.25} \qquad (14)$$

where $S_0$ is the solar constant at 1 AU, $\sigma$ is the Stefan-Boltzmann constant, $\alpha_s$ is the surface absorptivity in the solar spectrum, $\varepsilon_s$ is the surface emissivity in the infrared spectrum, and $\kappa$ is a coefficient which takes into account the surface area radiating in the infrared spectrum with respect to the one that receives the solar input. The maximum power that can be handled by the PPU (Power Processing Unit) is assumed to represent the upper limit for the engine thrust.

$$P_{in} = \min(P^*_{in}, P_{max}) \qquad (15)$$

The required power is a key factor in the sizing of the solar arrays and power system, and therefore can be used to provide an estimation of the overall dry mass of the spacecraft. The characteristics of the power system are summarised in Table 1.

**Table 1. Power system characteristics.**

| PARAMETER | $\eta_e$ | $\eta_S$ | $C_T$ | $T_0$ | $\kappa$ | $\varepsilon_s$ | $\alpha_s$ | $T_{max}$ | $P_{SS}$ |
|---|---|---|---|---|---|---|---|---|---|
| Value | 0.9 | 0.9 | 3e-4 K$^{-1}$ | 290 K | 1.8 | 1.0 | 0.8 | 423 K | 300 W |

## 3 - OPTIMISATION APPROACH

The design of an optimal trajectory involving a huge number of swing-bys with multiple resonant orbits was approached in two steps. First, we developed a procedure to find a set of good first guess solutions (FGS) both for the cruise part and for the capture part of the transfer. Afterwards, the most promising first guesses were optimised with a direct optimisation method. The trajectory model described in the previous section was used in the second step of the design while a simplified model was used in the first step of the design. In this section, we describe the procedure and the model employed to look for first guess solutions.

### 3.1 - GLOBAL SEARCH

The transfer from the Earth to Europa can follow different paths; each one characterised by a particular sequence of gravity assist manoeuvres, transfer time, launch date and departure $C_3$ (square of the asymptotic velocity at departure). In order to make the exploration of the space of all the possible transfers to Europa efficient, a reduced trajectory model was developed. The reduced model is based on a number of simplifying assumptions on the expected characteristics of the desired transfers. In particular, the orbits of the planets and of the moons were assumed to be coplanar and circular. Two search procedures were developed: one for the design of the transfer from the Earth to Jupiter, and one for the capture in the Jovian system and the tour of Jupiter's moons.

### 3.1.1 - *Interplanetary Transfers*

A key element of the design of a multi-gravity assist trajectory, either with low-thrust or chemical propulsion, is the position, along the orbit of the planets, at which the gravity manouvre occurs. The orbital position is a function of the time at which the manouvre occurs, therefore, we can build a multi-gravity assist trajectory by knowing the dates of the swing-bys. Given a sequence of swing-bys we want to find the set of dates at which we can perfom a gravity manouvre that minimises the Δv we have to deliver with the propulsion system. We call this problem *the phasing problem* in the remainder of the paper.

Solving exactly the phasing problem is a difficult task, therefore we try at first to compute an approximated solution to this problem. The solution is then used as first guess to design an optimal multi-gravity assist trajectory. First, we assume that the orbits of the planets are coplanar and that two subsequent gravity assist manoeuvres at planets A and B are connected by a conic arc. The exact time of transfer from A to B is approximated with a multiple of the Hohmann transfer time:

$$\Delta T = (2m+1)\Delta T_H = (2m+1)\pi\sqrt{(a_A + a_B)^3 / (8\mu_S)} \quad (16)$$

where $a_A$ and $a_B$ are the semi-major axes of planets A and B, and $m \in \mathbb{N}$ is the number of spacecraft revolutions on the transfer orbits. The assumption is that the exact time of the transfer will be in the interval $[\Delta T(1-\varepsilon) \; \Delta T(1+\varepsilon)]$ with $\varepsilon$ small. We compute now the departure times from planet A which solve the phasing problem.

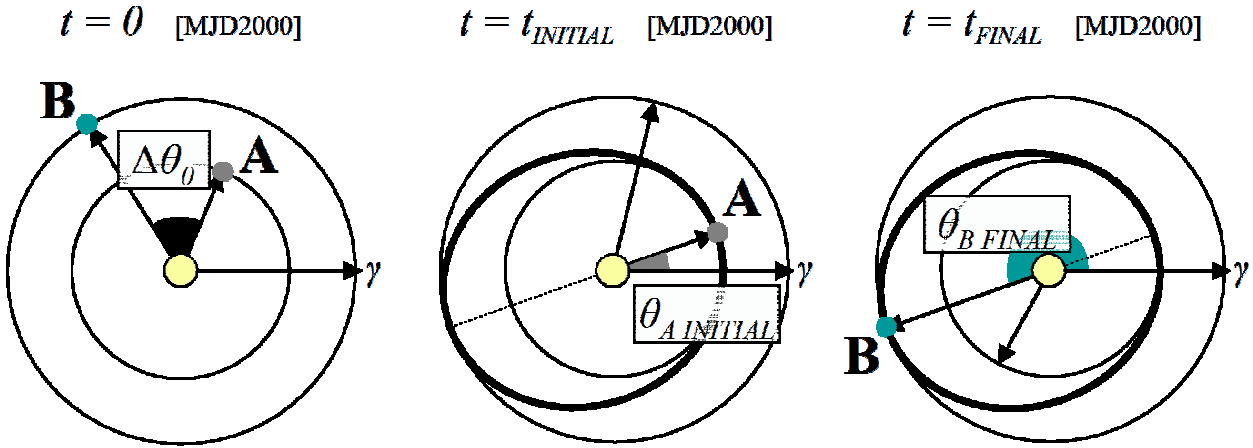

Figure 2. The phasing problem for a transfer from body A to body B.

We start considering Figure 2, which shows the position of the planets A,B at the reference time $t_0$ and at the initial and final time of the transfer from A to B. The angular position of the planet A,B at the initial and final time respectively are:

$$\begin{aligned}\theta_{AINITIAL} &= \theta_{A0} + \omega_A t_{INITIAL} \\ \theta_{BFINAL} &= \theta_{B0} + \omega_B t_{FINAL} = \theta_{B0} + \omega_B \Delta T + \omega_B t_{INITIAL}\end{aligned} \quad (17)$$

where $\omega_B$ and $\omega_A$ are the angular velocities of the arrival and departure planets respectively, and $\theta_{A0}$, $\theta_{B0}$ are their angular position at the reference time $t=0$. Now we impose the phasing contraint:

$$\theta_{BFINAL} - \theta_{AINITIAL} = (2k+1)\pi \qquad k \in \mathbb{N} \qquad (18)$$

Using the definitions of $\theta_{BFINAL}$ and $\theta_{AINITIAL}$ we find an expression for all the $t_{INITIAL}$ that solve the phasing problem:

$$t_{INITIAL}(k) = \frac{(2k+1)\pi - \omega_B \Delta T - \Delta\theta_0}{\omega_B - \omega_A} \qquad k \in \mathbb{N} \qquad (19)$$

The dates of departure can be further restricted so that the transfer occurs around the projection of the line of nodes onto the ecliptic.

Eq (19) solve the phasing problem providing the times $t_{INITIAL}$ at which the spacecraft leaves the planet A. We use now this equation to compute the launch date and swing-by dates for multiple gravity assist trajectories.

We start by dividing the trajectory in $N$ phases or transfer legs connecting subsequent swing-bys. For each phase $i$ we use equation (16) to compute the transfer time $\Delta T^{(i)}$ and Eq. (19) to compute the time $t^{(i)}_{INITIAL}$ at which the spacecraft is expected to leave the $i^{th}$ planet. To each $t^{(i)}_{INITIAL}$ we associate a departure date from Earth by subtracting the time of flight of the preceding phases.

$$t^{(i)}_{Launch}(k^{(i)}) = t^{(i)}_{INITIAL} - \sum_{j=1}^{i-1} \Delta T^{(j)} \qquad (20)$$

In an ideal scenario we would like to find a set of parameters $k^{(i)}$ so that all the $t^{(i)}_{Launch}$ are equal to the desired launch $t^{(1)}_{Launch}$. In general, however, we minimise the merit function:

$$F(k^{(1)},...,k^{(N)}) = \sum_{i=1}^{N} \left(t^{(i)}_{Launch} - t^{(1)}_{Launch}\right)^2 \qquad (21)$$

For each sequence of planets and set of spacecraft revolutions, problem (21) is solved, with a fast branch and prune procedure exploring the space of $N$ admissible integer variables $k^{(i)}$, $i = 1,...,N$. The resulting set of transfer solutions can be ordered according to the value of the merit function (21). The algorithm minimising the function $F$ in (21) was implemented in a sofware code called BS1.

### 3.1.2 - *Multiple Synchronous and Resonant Orbit Transfer*

The second part of the mission is a tour of the Galilean moons, which starts with a highly eccentric orbit around Jupiter (following the capture manoeuvre) and ends with a gravity assist at Europa. The trajectories in the Jupiter system use several Synchronous Orbit Tours (SOT) of the moons. A SOT is a sequence of resonant orbits linked by a number of gravity assist manoeuvres. Synchronous Orbit Tours are ballistic, there are no deep space manoeuvres and have no phasing constraints. Since a non-optimal sequence of resonant orbits (and gravity manoeuvres) can lead to very large transfer times, we developed the software tool BS2 to design minimum-transfer time SOTs.

### *3.1.2.1 - Synchronous Orbit tour*

A Synchronous Orbit Tour is a sequence of $N$ gravity assist manoeuvres with a given moon M at a constant position $\mathbf{r}_M$ along the moon's orbit. If we assume a linked-conic model for the $i^{th}$ gravity assist manoeuvre (with $i=1,...,N$), then we can define an incoming velocity vector $\mathbf{v}_i^-$ and an outgoing velocity vector $\mathbf{v}_i^+$ with respect to Jupiter at $\mathbf{r}_M$, before and after the GA. Once $\mathbf{v}_i^-$ and $\mathbf{v}_i^+$ are defined, we can compute the velocity $\tilde{\mathbf{v}}_i^-$ and $\tilde{\mathbf{v}}_i^+$ relative to the moon and the associated orbital plane of the gravity assist hyperbola. Since the physical constraint (8) must hold, in the remainder of this section we will use $\tilde{v}_i$ to indicate the modulus of the relative velocity without distinction between incoming or outgoing vectors. The initial conditions for a SOT are given by the vector $[\mathbf{r_M}, \mathbf{v}_1^-]^T$. The SOT terminates when the velocity modulus $v_N^+$ of the spacecraft reaches a target value $v_{Target}$ (or equivalently when the orbital energy of the spacecraft reaches a target value). The outgoing velocity $\mathbf{v}_i^+$ uniquely defines the orbit connecting the gravity assist number $i$ to the gravity assist number $i+1$. Each orbit has a unique period $T_i$ and a unique resonance ratio $\rho_i = n_i/m_i = T_i/T_M$, where $T_M$ is the period of the moon and $n_i$, $m_i$ are two integer numbers. The transfer time associated to each orbit is $\Delta t_i = n_i T_M = m_i T_i$, with the semi-major axis equal to:

$$a_i = \rho_i^{2/3} a_M \quad (22)$$

where $a_M$ is the semi-major axis of the orbit of the moon. Note that, we assume the orbits of the moons to be circular therefore $r_M = a_M$.

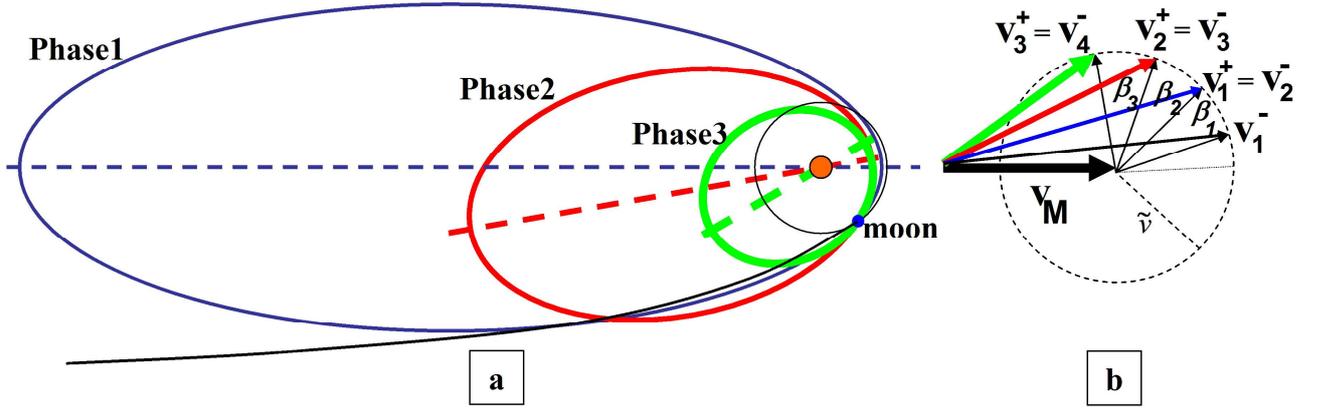

Figure 3. Schematic of the SOT.

Figure 3a shows the first phase of a synchronous orbit tour. Figure 3b shows that the relative velocity $\tilde{\mathbf{v}}_i^-$ is always in the plane containing the initial spacecraft velocity $\mathbf{v}_1^-$ and the velocity of the moon $\mathbf{v}_M$. Note that as $\tilde{\mathbf{v}}_i^+$ rotates at each GA, the velocity $\mathbf{v}_i^+$ decreases. For a synchronous tour, the outgoing velocity $\mathbf{v}_i^+$ is subject to the following constraints:

$$v_i^+ = \sqrt{\frac{\mu}{a_M}}\sqrt{2 - \rho_i^{-2/3}} \quad (23)$$

$$v_i^+ = \sqrt{v_M^2 + \tilde{v}_i^2 + 2v_M \tilde{v}_i \cos\left(\alpha_0 + \sum_{j=1}^{i} \beta_j\right)} \quad (24)$$

where the deviation angle $\beta_i$ is bounded by:

$$\beta_0 \leq \beta_i \leq \beta_{MAX} \quad (25)$$

and $\beta_{MAX}$ is given by Eq. (10) with the minimum altitude of the GA. Furthermore, if we define:

$$v_{\min,i} = v_i^+\big|_{\beta_i=\beta_{MAX}} \quad (26)$$

$$v_{MAX,i} = v_i^+\big|_{\beta_i=0} = v_{i-1} \quad (27)$$

then the modulus of the velocity vector is bounded by:

$$v_{\min,i} \leq v_i^+ \leq v_{MAX,i} \quad (28)$$

The combination of Eqs.(23)-(28) yields:

$$\rho_{\min,i} \leq \rho_i \leq \rho_{MAX,i} \tag{29}$$

with

$$\rho_{MAX,i} = \frac{n_{i-1}}{m_{i-1}} \tag{30}$$

and

$$\rho_{\min,i} = f(v_i) = \left\{ 2 - \frac{a_M}{\mu} \left[ v_M^2 + \tilde{v}_i^2 + 2v_M \tilde{v}_i \cos\left( \beta_{\max} + \mathrm{a}\cos\left( \frac{\frac{\mu}{a_M}\left(2 - \rho_{i-1}^{-2/3}\right) - v_M^2 - \tilde{v}_i^2}{2 v_M \tilde{v}_i} \right) \right) \right] \right\}^{-3/2} \tag{31}$$

Equation (31) can be satisfied only when $m_i > m_{MIN,i}$, where $m_{MIN,i}$ is defined as:

$$m_{MIN,i} = \min_{I_m \neq \varnothing} m_i \tag{32}$$

and the set $I_m$ is defined as:

$$I_m = \left\{ j \in \mathbb{N} \,\middle|\, m_i \rho_{\min,i} \leq j \leq m_i \rho_{MAX,i} \right\} \tag{33}$$

Also, $n_i \geq n_{MIN,i}$, where $n_{MIN,i}$ is the smallest integer greater than $m_i \rho_{\min,i}$. As a special case, the velocity after the last GA, $v_N$, only needs to satisfy the constraint $v_N = v_{Target}$. Thus, the sequence ends when $v_{\min,i} < v_{Target}$, or equivalently $\rho_{\min,i} > f(v_{Target})$. Finally, the transfer time of the entire SOT is:

$$T = \sum_{i=1}^{N} \Delta t_i = T_M \sum_{i=1}^{N} n_i \tag{34}$$

Now an optimal SOT is such that, given the initial conditions $\left[\mathbf{r_M}, \mathbf{v}_1^-\right]^T$, the time to achieve the target velocity is minimal. Equivalently, looking at Eq. (34), an optimal SOT is such that:

$$\min_{\mathbf{n},\mathbf{m},N \in \mathbb{N}} \sum_{i=1}^{N} n_i \tag{35}$$

subject to

$$\rho_{\min,i} < \frac{n_i}{m_i} < \frac{n_{i+1}}{m_{i+1}} \tag{36}$$

$$m_i \geq m_{MIN,i} \tag{37}$$

$$m_i \geq m_{MIN,i} \tag{38}$$

$$n_i \geq n_{MIN,i} \tag{39}$$

$$\rho_{\min,N} > f(v_{Target}) \tag{40}$$

Eqs. (35)-(40) defines an integer programming problem which has a number of solutions that grown exponentially with the number of gravity assists. In order to reduce the solution space we seek for solutions that maximise the GA efficiency:

$$\eta_i = \frac{\beta_i}{\beta_{MAX}} \tag{41}$$

The GA efficiency $\eta_i$ of the phase $i$ increases if the angle $\beta_i$ increases. Eq. (24) and Figure 3b show that $\beta_i$ increases when $\tilde{v}_i^+$ increases, thus when the resonant ratio $\rho_i$ decreases (Eq. 23). We can conclude that the GA efficiency $\eta_i$ increases with either an increase of $m_i$ or a decrease of $n_i$:

$$\begin{cases} \dfrac{\Delta \eta_i}{\Delta n_i} < 0 \\ \dfrac{\Delta \eta_i}{\Delta m_i} > 0 \end{cases} \quad (42)$$

Since choosing a low $n_i$ implies also choosing a lower transfer time, we replace the inequality constraint $n_i \geq n_{MIN,i}(m_i)$ with:

$$n_i = n_{MIN,i}(m_i) \quad (43)$$

With this assumption, the solution space is significantly reduced and can be explored using a branch and bound technique. In particular, the software tool BS2 implements the search, and prunes a branch when $n_i \geq n^{OPT} - \sum_1^{i-1} n_j$ where $n^{OPT}$ is the current optimum solution. BS2 was implemented in Matlab and finds a six-gravity-assist SOT in less than a second on a 2GHz laptop.

### 3.1.2.2 - The transfer to the next moon

In the previous section we showed how BS2 computes the fastest SOT to reach a target velocity $v_{Target}$. Often the main purpose of the SOT is to reduce the energy in order to minimize the relative velocity of the spacecraft with respect to another moon. Figure 4 shows the geometry of a transfer from moon A to moon B. The aim is to minimize the velocity of the spacecraft relative to a Moon ($\tilde{v}_B$ in Figure 4) subject to the constraints:

$$\|\tilde{v}_A\| = \text{fixed} \quad (44)$$

$$\alpha_{min} < \alpha < \pi \quad (45)$$

where $\alpha_{min}$ is the angle corresponding to a transfer from A to B tangent to the arrival moon B. Note that, when $\alpha < \alpha_{min}$ the transfer trajectory does not intersect the orbit of the moon B.

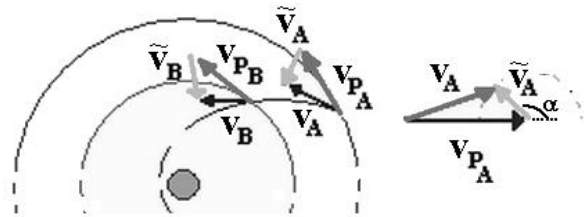

Figure 4. Schematic of a trajectory intersecting two orbits.

The problem was solved analytically as follows. First we recall that the relative velocities at the moons A and B are:

$$\tilde{v}_B^2 = v_{M_B}^2 + v_B^2 - 2v_B^\theta v_{M_B} \quad (46)$$

$$\tilde{v}_A^2 = v_{M_A}^2 + v_A^2 - 2v_A^\theta v_{M_A} \quad (47)$$

where the velocity of a moon M on its orbit is $v_M = \sqrt{\mu/r_M}$ and $v^\theta = v_M + \tilde{v}\cos\alpha$ is the transverse velocity of the spacecraft at the moon. Now if we subtract (47) from (46) and use the conservation of the energy and of the angular momentum, we get:

$$\tilde{v}_B^2 = \tilde{v}_A^2 + 3\mu\left(\frac{1}{r_A}-\frac{1}{r_B}\right) - 2v_{M_A}^2\left[\left(\frac{r_A}{r_B}\right)^{3/2}-1\right] - 2\tilde{v}_A^2 v_{M_A}\left[\left(\frac{r_A}{r_B}\right)^{3/2}-1\right]\cos\alpha = k_1 - k_2\cos\alpha \tag{48}$$

where the term $k_2$ is positive, therefore the velocity $\tilde{v}_B^2$ is minimal if $\alpha = \alpha_{\min}$, i.e. when the transfer is tangent to the moon B.

### 3.1.3 - *Optimal Impulsive Transfer Trajectories*

The algorithm implemented in BS1 is fast and give a good estimation of the set of optimal sequences of swing-by manoeuvres. In particular, it gives a good estimation of the right phasing of the encounters with the planets. However, it does not provide the correct computation of the total $\Delta v$. Therefore, a more sophisticated trajectory model, including deep-space manoeuvres, was implemented in a third algorithm, called BS3. BS3 was used to generate a more accurate FGS for all the electric propulsion options. The underlying idea is that optimal impulsive transfers can be regarded as a limit case of a minimum mass, low-thrust transfer with no limit on the thrust level. Therefore, we assume that the solution that minimises the total impulsive $\Delta v$ will also minimise the duration of the low-thrust arcs and therefore the corresponding propellant consumption. As the results of our case study demonstrate, this assumption is true, provided that there is enough time to spread the velocity variation over the whole arc with the given level of thrust. A similar approach, was recently used by Okutsu et al. [9] to design optimal low-thrust transfers to Jupiter, with analogous results.

BS3 uses three dimensional analytical ephemerides for the positions of the planets and a linked conic approximation for gravity assist manoeuvres. The trajectory is split into phases connecting two consecutive GAs, and each phase is again split into two subarcs connected by a deep space maneuvers (DSM). Figure 4. shows the geometry of an arc connecting the gravity assist $i$ at planet $P_i$ to the gravity assist $i+1$ with planet $P_{i+1}$. The independent variables are the time of the two gravity assists $t_i$ and $t_{i+1}$ and the time $t_{DSM_i}$ and position $\mathbf{r}_{DSM_i}$ of the deep space maneuver. Each subarc is then computed as the solution of the Lambert's problem[14]. The vectors. $\mathbf{r}_{P_i}$ and $\mathbf{v}_{P_i}$ are the position and velocity of the planet $P_i$ while the vectors $\mathbf{v}^-_{DSM_i}$ and $\mathbf{v}^+_{DSM_i}$ are the velocities of the spacecraft before and after the DSM.

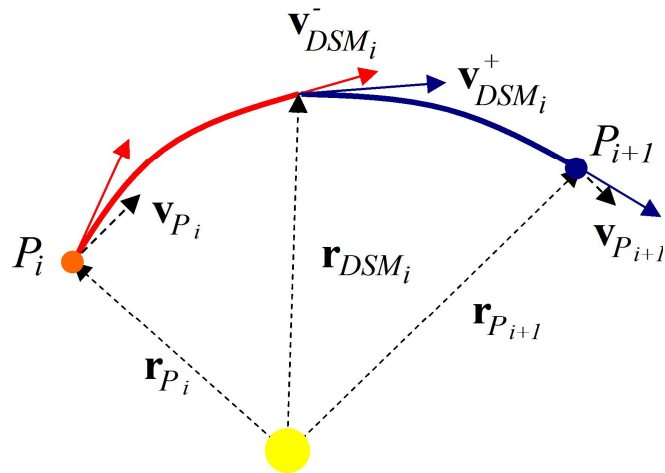

Figure 5. Schematic of the a single trajectory leg implemented in BS3.

Once all the subarcs are assembled together, the problem is to find the date of each planet encounter, the launch and arrival dates, and the points in space and time where to perform a manoeuvre in order to minimise the sum of all required $\Delta vs$. The trajectory that minimises the fuel consumption from planet $P_1$ to $P_N$ makes use of (*N*-1) impulsive manoeuvres and (*N*-2) flybys, therefore the solution vector of the resulting optimisation problem has (6*N*-6) components:

$$\mathbf{X} = [t_{P_1},...,t_{P_N},\ t_{DSM_1},...,t_{DSM_{N-1}},\ r_{DSM_1},...,r_{DSM_{N-1}},\ \tilde{r}_{p,1},\ ...,\ \tilde{r}_{p,N-2}]^T \tag{49}$$

where $t_{Pi}$ is the time at the $i^{th}$ planetary encounter, $t_{DSMi}$ is the time of the $i^{th}$ impulsive manoeuvre, $r_{DSMi}$ is the position vector of the $i^{th}$ impulsive manoeuvre and $\tilde{r}_{p,i}$ is the pericentre altitude of $i^{th}$ GA manoeuvre.

The optimisation problem can be formulated as follows

$$\min_{\mathbf{X} \in D \subseteq \Re^{(6N-6)}} F_{BS3}(\mathbf{X}) = \sum_{i=1}^{N-1} \Delta v_i = \sum_{k=1}^{N-1} \left\| \mathbf{v}^+_{DSM_i} - \mathbf{v}^-_{DSM_i} \right\| \tag{50}$$

subject to constraints (8) and (9). Two additional constraints were introduced in order to limit the departure velocity from Earth and the relative arrival velocity at Jupiter.

$$\tilde{v}_1 < \tilde{v}_{departure} \qquad \tilde{v}_N < \tilde{v}_{arrival} \tag{51}$$

From BS1 we obtained a first guess value for the time components in the vector **X**, then an NLP solver was used to find an optimal solution to problem (50).

### 3.2 - DIRECT OPTIMISATION

Once a promising first guess solution, both for the cruise part and for the Jovian tour part, was available, the corresponding trajectory was transcribed with a direct approach based on Finite Elements in Time. This approach decomposes a general trajectory design problem into $M$ phases, each one characterised by a time domain $D^j$ with $j=1,...,M$, a set of $m$ dynamic variables **x**, a set of $n$ control variables **u** and a set of $l$ parameters **p**. Each phase $j$ has an objective function:

$$J^j = \phi^j\left(\mathbf{x}_0^b, \mathbf{x}_f^b, t_f, \mathbf{p}\right) + \int_{t_i}^{t_f} L^j(\mathbf{x}, \mathbf{u}, \mathbf{p}) dt \tag{52}$$

a set of dynamic equations:

$$\dot{\mathbf{x}} - \mathbf{F}^j(\mathbf{x}, \mathbf{u}, \mathbf{p}, t) = 0 \tag{53}$$

a set of algebraic constraints on states and controls:

$$\mathbf{G}^j(\mathbf{x}, \mathbf{u}, \mathbf{p}, t) \geq \mathbf{0} \tag{54}$$

and a set of boundary constraints:

$$\psi^j\left(\mathbf{x}_0^b, \mathbf{x}_f^b, \mathbf{p}, t\right)\Big|_{t_0}^{t_f} \geq 0 \tag{55}$$

A special subset of boundary constraints, called inter-phase link constraints, are used to assemble all the phases together:

$$\psi^j\left(\mathbf{x}_j^b, \mathbf{x}_{j-1}^b, \mathbf{p}, t\right) \geq 0 \tag{56}$$

The time domain $D(t_0, t_f) \subset \Re$ relative to each phase $j$ is further decomposed into $N$ finite time elements $D^j = \bigcup_{i=1}^N D_i^j(t_{i-1}, t_i)$. For each time element $D^j_i$, states and controls [**x**, **u**] are parameterised as follows:

$$\begin{Bmatrix} \mathbf{x} \\ \mathbf{u} \end{Bmatrix} = \sum_{s=1}^{p} f_s(t) \begin{Bmatrix} \mathbf{x}_s \\ \mathbf{u}_s \end{Bmatrix} \tag{57}$$

where the basis functions $f_s$ are chosen within the polynomial space of order $(p-1)$:

$$f_s \in P^{p-1}(D_i^j) \qquad (58)$$

Therefore we can define a finite element as composed of a sub-domain $D^j_i$ and a subset of parameters $[\mathbf{x}_s, \mathbf{u}_s, \mathbf{p}]$. A group of finite elements forms a phase, and a group of phases forms the original trajectory. Note that an additional parameter set $\mathbf{p}$ may occur in all constraint equations depending on their function in the optimisation problem. Furthermore, phases can be put in sequential order or in parallel with other phases depending on their time domain and inter-phase link constraints. Thus, two phases can share the same time domain but have different parameterisations, different dynamic model and different objective functions.

Now, considering a general phase $j$, in order to integrate the differential constraints (53), on each finite element $i$, differential equations are transcribed into a weighted residual form considering boundary conditions of the weak type:

$$\int_{t_i}^{t_{i+1}} \left\{ \dot{\mathbf{w}}^T \mathbf{x} + \mathbf{w}^T \mathbf{F}^j \right\} dt - \mathbf{w}_{i+1}^T \mathbf{x}_{i+1}^b + \mathbf{w}_i^T \mathbf{x}_i^b = 0 \qquad i = 1, \ldots, N-1 \qquad (59)$$

where $\mathbf{w}(t)$ are the generalised weight (or test) functions defined as:

$$\mathbf{w} = \sum_{s=1}^{p+1} g_s(t) \mathbf{w}_s \qquad (60)$$

and $g_s$ is taken within the polynomial space of order $p$:

$$g_s \in P^p(D_i^j) \qquad (61)$$

Now the problem is to find the vectors $\mathbf{x}_s \in \Re^{p*m}$, $\mathbf{u}_s \in \Re^{p*n}$, $\mathbf{p} \in \Re^l$ and $[\mathbf{x}^b_f, \mathbf{x}^b_0] \in \Re^m$ that satisfy the variational equation (59) along with algebraic and boundary constraints:

$$\mathbf{G}^j(\mathbf{x}, \mathbf{u}, \mathbf{p}, t) \geq 0 \qquad (62)$$

$$\psi^j\left(\mathbf{x}_0^b, \mathbf{x}_f^b, \mathbf{p}, t\right)\bigg|_{t_0}^{t_f} = 0 \qquad (63)$$

where the quantities $\mathbf{x}_s$ and $\mathbf{u}_s$ are internal node values, $\mathbf{x}^b_f$ and $\mathbf{x}^b_0$ are the boundary values. Note that generally, the order $p$ of the polynomials can be different for states and controls.

Each integral of the continuous forms (52) and (59) is then replaced by a $q$-point Gauss quadrature sum, where $q$ is taken equal to $p$. For a continuous solution, in order to preserve the continuity of the states at matching points, the following condition must hold:

$$\mathbf{x}_i^b = \mathbf{x}_{i+1}^b \qquad i = 1, \ldots, N-2 \qquad (64)$$

Thus all the boundary quantities in Eq. (64) cancel one another except for those at the initial and final times. Algebraic constraint equation (62) can be collocated directly at the Gauss nodal points:

$$\mathbf{G}^j_s(\mathbf{x}_s(\xi_s), \mathbf{u}_s(\xi_s), \mathbf{p}, \xi_s) \geq 0 \qquad (65)$$

The resulting set of non-linear algebraic equations, assembling all the phases, along with discretised objective function (52) can be seen as a general non-linear programming problem (NLP) of the form:

$$\min J(\mathbf{y}) \qquad (66)$$

subject to:

$$\begin{aligned} \mathbf{c}(\mathbf{y}) &\geq 0 \\ \mathbf{b}_l &\leq \mathbf{y} \leq \mathbf{b}_u \end{aligned} \qquad (67)$$

where $\mathbf{y} = \left[\mathbf{x}_s, \mathbf{u}_s, \mathbf{x}_0^b, \mathbf{x}_f^b, t_0, t_f, \mathbf{p}\right]^T$ is the vector of NLP variables, $J(\mathbf{y})$ is the objective function to be minimised, $\mathbf{c}(\mathbf{y})$ a vector of non-linear constraints and $\mathbf{b}_l$ and $\mathbf{b}_u$ are, respectively, the lower and upper bounds on NLP variables. The DFET approach was implemented in a software code called DITAN (Direct Interplanetary Trajectory ANalysis)[3][4][15]=. Note that the DFET formulation accommodates the optimisation of the controls parameters and a set of static parameters $\mathbf{p}$ that can be related to spacecraft system characteristics (e.g. $I_{sp}$, power level, etc.).

## 4 - MISSION DESIGN

In this section, we present some results that illustrate how to apply the numerical approaches, proposed in the previous section, to the design of a low-thrust multigravity assist trajectory to Europa. The tools BS1, BS2, BS3 and DITAN will be applied according to the following procedure:

1. BS1 is run to generate a pool of candidate solutions that present a minimum violation if the phasing constraint (18). The most promising ones are optimised with BS3 imposing a limit on the $C_3$ at launch, the arrival at the sphere of influence of Jupiter but no constraint on the arrival velocity. The solutions computed with BS3 are passed to DITAN to verify the assumption on the validity of impulsive solutions as first guesses for low-thrust solutions. The results of this first part of the design are illustrated in section 4.2.

2. BS2 is run to find an optimal sequence of resonant flybys to reach Europa (see section 4.3). Prior to this, we performed an analysis of the required conditions for a capture in the Jupiter system (see section 4.3.1).

3. The solutions computed with BS3 for the Earth-Jupiter leg, are then optimised with DITAN imposing an arrival velocity at the sphere of influence compatible with the results of BS2.

4. The resonant orbits designed with BS2 are added, in DITAN, to the rest of the Earth-Jupiter transfer, and the trajectory is optimised as a whole (see section 4.4).

### 4.1 - SYSTEM PARAMETER DEFINITION

The dynamic model in Eq. (1) includes the initial mass of the spacecraft and the maximum thrust level that the engines can deliver. In particular, the required thrust level at Jupiter is the driving parameter for the selection of the engine. We considered a minimum of 35 mN of thrust at Jupiter in order to be captured (see the following analysis on the capture manoeuvre) and a thrust-to-power ratio of 27 mN/W for the engine. The power required at end of life (EOL) at Jupiter was estimated using the power-to-thrust ratio of the engine and Eq.(13). From the EOL power, the beginning of life power (BOL) at 1 AU was used to calculate the size of the solar arrays and their mass. We considered triple junction GaAs cells combined with solar concentrators in order to reduce panel size. We assumed a total mass of the spacecraft, including propellant, of 1500 kg and a maximum $C_3 = 3.16$ km$^2$/s$^2$. The final configuration consisted of seven 150 mN ion thrusters for a total mass of the power and propulsion systems reported in Table 2.

**Table 2. Power and Propulsion System Budget.**

| Required Power | 40 kW (EOL) |
|---|---|
| Solar Array Area | 133 m$^2$ |
| Ion Thruster | 49 kg |
| PPU | 97 kg |
| Tank+Harness+Piping | 30 kg |
| Flux Lines | 15 kg |

| Control Unit | 5 kg |
| --- | --- |
| Propellant | 187 kg |
| Solar Array | 416 kg |
| Total SEP Mass | 773 kg |
| Total Spacecraft Mass | 1500 kg |

**4.2 - EARTH TO JUPITER TRANSFER**

At first, we investigated the transfers from the Earth to Jupiter. Venus was chosen as the first planet, in the sequence of gravity assist manoeuvres to reach Jupiter, since the model in BS1 does not deal with Earth-Δv-Earth arcs, therefore an Earth-Earth resonant flyby could not be used. On the other hand, Venus provides a higher Δv change with respect to the other remaining inner planets. Afterwards a flyby strategy to raise the apocentre and pericentre was devised based on the algorithm coded in BS1. BS1 yielded a wide range of transfers with several sequences. However, only the two most promising solutions generated by BS1 were optimised with the model implemented in BS3 and then passed to DITAN. The requirement for the optimised solution was to arrive at the sphere of influence of Jupiter with no limit on the arrival velocity. In particular we considered the sequence EVMEJ, and BS1 found solutions which included a phasing orbit before the encounter with Venus. The results are shown in Figure 6 and Figure 7 and Table 3 and Table 4. The left parts of Figure 6 and Figure 7 represent the first guess solution computed with BS1 while the right hand sides of the two figures represent the optimised solutions with DITAN. The figures show that the optimised solution is very close to the first guess, with very limited thrust arcs. The usefulness of the search performed by BS1 is illustrated also by Table 3 and Table 4 where the dates of arrival and departure found by BS1 are compared to the optimised dates of encounter with the planets computed with DITAN. Considering that DITAN uses the JPL DE405 ephemeris while BS1 uses the mean orbital elements, the first guess and the optimised solution are very similar except for the arrival at Jupiter of the first solution.

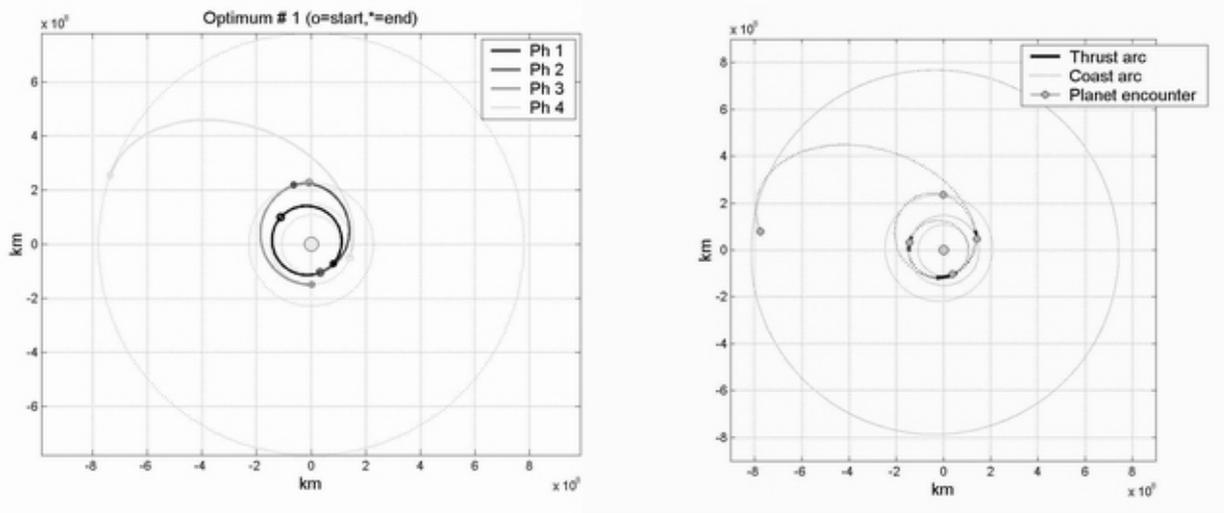

Figure 6. Solution number 1: First Guess Solution generated by BS1 (left) and optimised with electric propulsion (right).

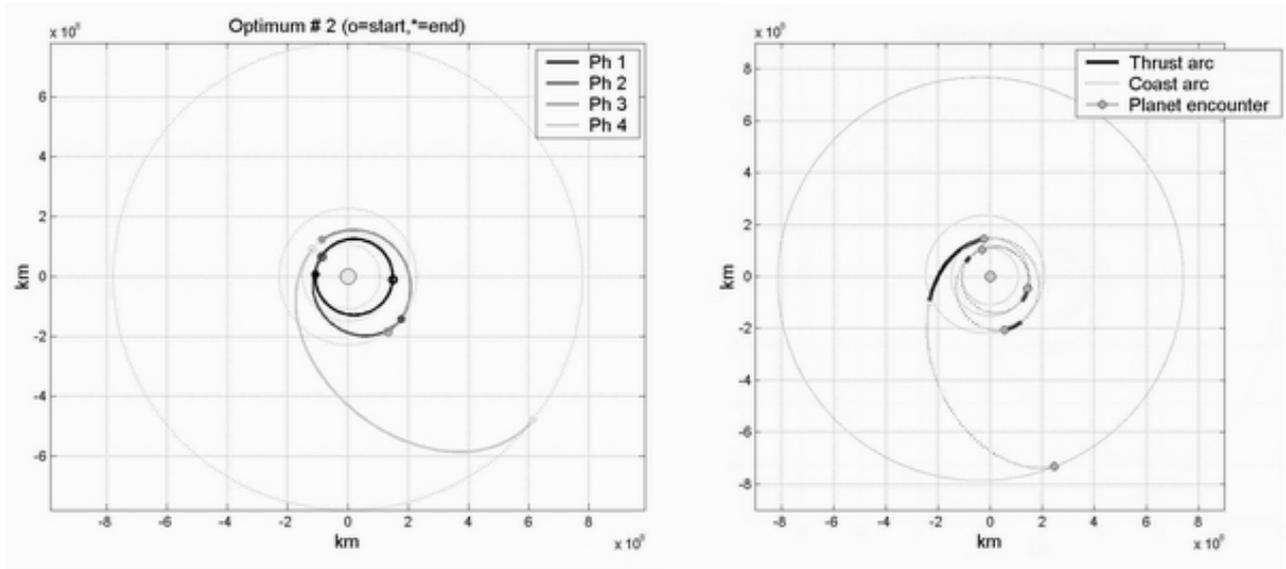

Figure 7. Solution number 2 : FGS generated by BS1 (left) and optimised with electric propulsion (right).

Table 3. Optimal solution 1.

| GA PLANET | FIRST GUESS | OPTIMISED |
|---|---|---|
| Earth | Dep (Phase 1): 12-02-2010 (3692 MJD2000) | Dep:07-03-2010 (3719 MJD2000) |
| Venus | Arr (Phase 1): 22-04-2011 (4130 MJD2000) Dep (Phase 2):02-04-2011 (4110 MJD2000) | Flyby: 04-04-2011 (4112 MJD2000) |
| Mars | Arr (Phase 2): 06-11-2011 (4328 MJD2000) Dep (Phase 3): 08-10-2011 (4299 MJD2000) | Flyby: 23-09-2011 (4284 MJD2000) |
| Earth | Arr (Phase 3): 23-06-2012 (4558 MJD2000) Dep (Phase 4): 01-09-2012 (4628 MJD2000) | Flyby: 10-10-2012 (4667 MJD2000) |
| Jupiter | Arr (Phase 4): 26-05-2015 (5625 MJD2000) | Arr: 27-04-2016 (5962 MJD2000) |

Table 4. Optimal solution 2.

| GA PLANET | FIRST GUESS | OPTIMISED |
|---|---|---|
| Earth | Dep (Phase 1): 17-09-2003 (1356 MJD2000) | Dep: 04-09-2003 (1343 MJD2000) |
| Venus | Arr (Phase 1): 28-11-2004 (1794 MJD2000) Dep (Phase 2): 07-11-2004 (1773 MJD2000) | Flyby: 16-10-2004 (1751 MJD2000) |
| Mars | Arr (Phase 2): 12-06-2005 (1990 MJD2000) Dep (Phase 3): 12-05-2005 (1959 MJD2000) | Flyby: 25-04-2005 (1942 MJD2000) |
| Earth | Arr (Phase 3): 26-01-2006 (2218 MJD2000) Dep (Phase 4): 11-02-2006 (2234 MJD2000) | Flyby: 29-12-2005 (2190 MJD2000) |
| Jupiter | Arr (Phase 4): 05-11-2008 (3232 MJD2000) | Arr: 24-07-2008 (3128 MJD2000) |

**4.3 - TOUR OF THE JOVIAN SYSTEM**

The second part of the transfer comprises a capture phase and a descent to Europa made of a number of SOTs.

**4.3.1 -** *Capture*

At arrival at Jupiter the spacecraft exploits a gravity assist manoeuvre with one of the moons to gain the required $\Delta v$ to be captured in the Jovian system. We compute the period of the post swing-by orbit, for each of the principal moons, as a function of the velocity at the sphere of influence $v_{SOI}$ and of the angle $\gamma$, between the incoming velocity $v^-$ and the velocity of the moon $\mathbf{v}_M$. The geometry is shown in Figure 8.

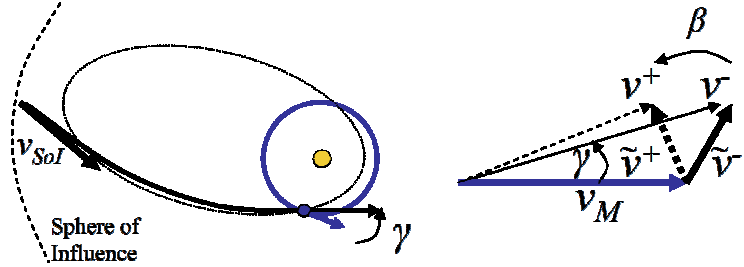

Figure 8. The first swing-by at a Galilean moon captures the spacecraft in the Jovian system. On the left, the hyperbolic approach and the captured orbit (dashed curve). On the right, geometry of the swing-by.

The period of the post swing-by orbit can be computed as a function of the orbital energy post swing-by $E$:

$$T(\gamma, v_{SoI}) = 2\pi\mu_J \left(-2E\right)^{-3/2} \tag{68}$$

with

$$E(\gamma, v_{SoI}) = \frac{1}{2}(v^-)^2 - v_M(v^-\cos\gamma + v_M \cos\beta) + v_M v^- \cos(\gamma + \beta) \tag{69}$$

and

$$\beta = 2\arcsin\left(\frac{\mu_M}{\mu_M + \tilde{r}_p\left((v^-)^2 + v_M^2 - 2v_M v^- \cos\gamma\right)}\right) \tag{70}$$

$$v^- = \sqrt{v_{SoI}^2 + 2v_M^2 - 2\frac{\mu_J}{r_{SoI}}} \tag{71}$$

where $\mu_J$ is the gravity constant of Jupiter, $\mu_M$ is the gravity constant of the moon, $v_M$ is the modulus of the velocity of the moon, $v^-$ is the modulus of the incoming velocity of the spacecraft, $\tilde{r}_p$ is the pericentre of the swing-by hyperbola. Figure 9 shows the contour plot of the periods of the post swing-by orbits for a pericentre altitude at each flyby of 300 km. Examining the period of the first orbit after the capture, Ganymede, appears to be the best planet for the first synchronous orbit. Figure 10 shows the limit curves for which the post swing-by orbit is parabolic. In particular, the figure shows that a Ganymede swing-by can insert the spacecraft into a parabolic orbit around Jupiter if the velocity at the sphere of influence is below 4.8 km/s(around 4.1 km/s at infinity). Also, Figures 9 and 10 show that the best capture conditions are achieved when γ is around 15°-20°, regardless of the velocity at the sphere of influence, and regardless of the moon chosen for the first swing-by.

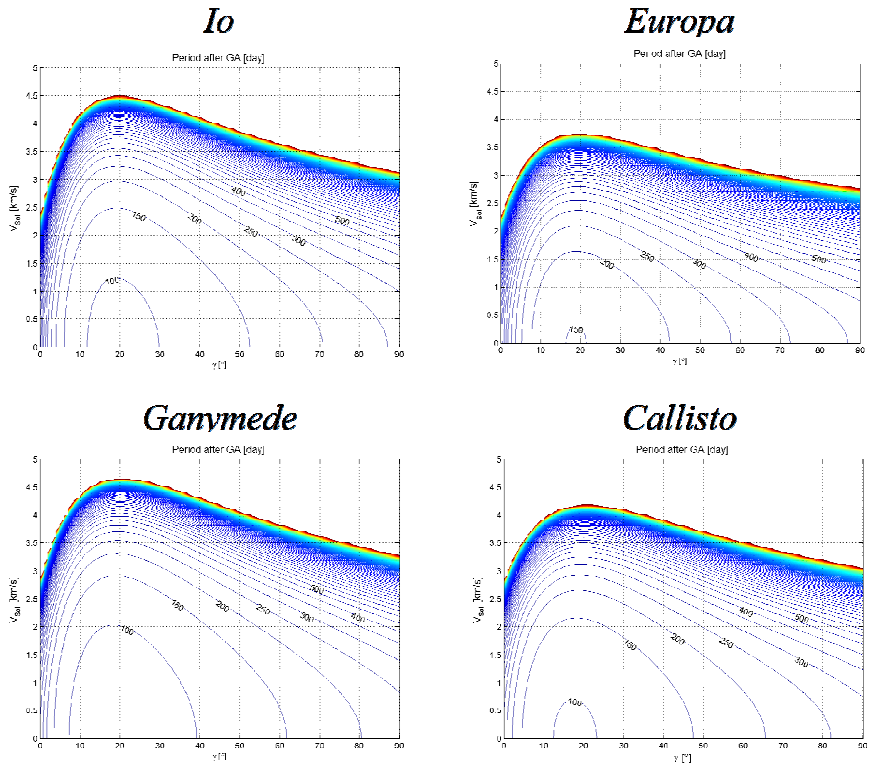

Figure 9. Possible period of the first orbit after the first gravity assist of one of the moons as a function of the velocity modulus at the sphere of influence of Jupiter and of the angle γ.

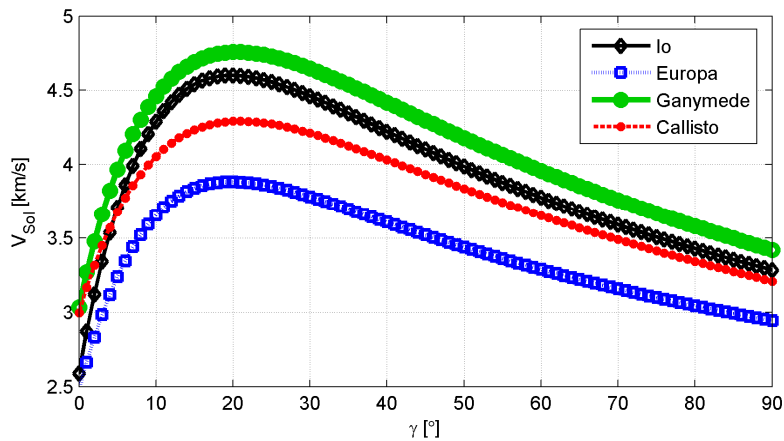

Figure 10. Maximum velocity at the sphere of influence to insert the spacecraft into a parabolic orbit with a moon swing-by, as a function of the angle γ for different moons.

#### 4.3.2 - *Resonant Descent to Europa*

If we consider an ideal planar bi-impulsive transfer between Ganymede and Europa, both having circular orbits, the lowest velocity relative to Europa is $\tilde{v}_{Eu} = 1.49$ km/s. In general, however, the relative velocity at Ganymede is higher then the one required for the Hohmann transfer. Yet low relative velocities at Europa (in the range from 1.5 km/s to 2 km/s) can be achieved by a synchronous tour leading to a nearly tangent-to-target transfer. This tour would exploit only gravity assist manoeuvres of the moon, minimising the use of the low-thrust engine. The moon sequence is made up of a synchronous Ganymede tour until the pericentre has dropped to the orbit of Europa. Then a synchronous Europa tour is performed until the apocentre has decreased down to the orbit of Ganymede. A final Ganymede flyby puts the spacecraft into a Ganymede-Europa tangent orbit.

The tour of the Jovian moons was designed using BS2 to estimate the best sequence of swing-bys. Once the initial conditions for the tour were derived, they were passed to DITAN as terminal constraints for the transfer trajectory from the Earth to Jupiter. The result is represented in Figure 11. The same sequence estimated using BS2 was then optimised using DFET. The result is reported in Table 5 for the tour of Ganymede and Table 6 for the tour of Europa. The last column represents the resonances after the optimisation with DITAN while the fourth column represents the resonances computed with BS2. For example, for the first phase BS2 predicts 75 revolutions of Ganymede (the value of *n*) and one of the spacecraft while the optimised solution is a bit faster with 69 revolutions of the moon and one of the spacecraft.

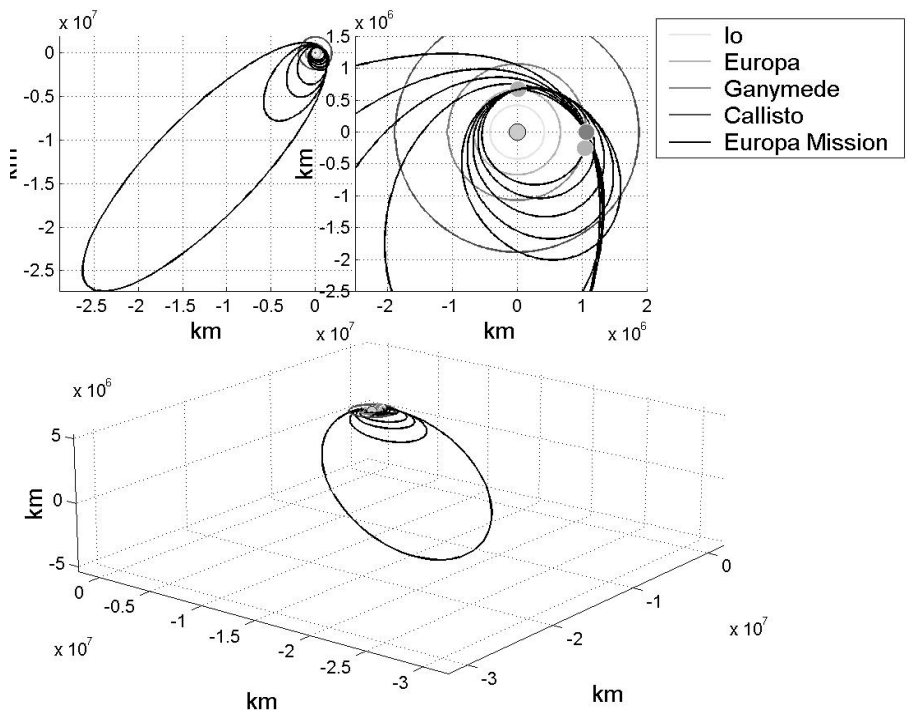

Figure 11. First Guess Solution for the tour of the Jovian moons.

Since the real motion of the planets differs slightly from the mean motion, there is an error of 0.15 % in the calculated distance from the spacecraft to Jupiter while approaching Ganymede. Such a small error is amplified when evaluating the semi-major axis. The error in *a* is 5-6 % for a high energy orbit and decreases to 0.5-0.6 % for lower energy orbits. This non-linear behaviour explains the discrepancies between the first guess and optimised solutions.

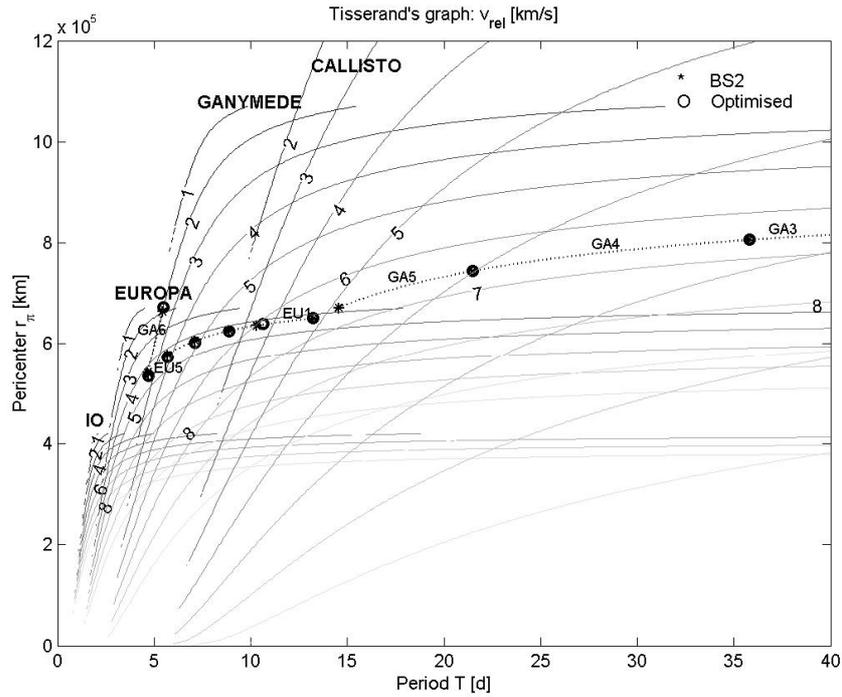

Figure 12. Tour representation on Tisserand's plane.

As a further confirmation of the effectiveness of the strategy implemented in BS2, the designed tour of Jovian moons was represented on Tisserand's plane[16] (see Figure 12), along with the optimised sequence. The two sequences are remarkably similar, which suggests that the model implemented in BS2 provides very good first guess solutions for the design of a synchronous tour of Jupiter's moons. It should be noted that the tour proposed here does not take into account any constraints on the radiation dose. However, by using BS2, the total dose can be easily minimised by imposing further pruning conditions. The distinct advantage of this sequence is that it allows the spacecraft to reach Europa with a relative velocity of 1.65 km/s (optimised solution). This is a low value comparable to the best results found in literature[17][18]. If required, endgame strategies can be applied at this point to further reduce the cost of the orbit insertion manoeuvre at Europa[17].

Table 5. Optimal Ganymede sequence from BS2.

| PHASE NUMBER | APOCENTER (FGS) | APOCENTER (OPTIMISED) | RESONANCE n:m (FGS) | RESONANCE n:m (OPTIMISED) |
|---|---|---|---|---|
| 1 | 1.90e7 km | 1.80e7 km | 75:1 | 69:1 |
| 2 | 4.97e6 km | 5.92e6 km | 10:1 | 7:1 |
| 3 | 3.13e6 km | 3.19e6 km | 5:1 | 5:1 |
| 4 | 2.22e6 km | 2.22e6 km | 3:1 | 3:1 |

Table 6. Optimal Europa sequence from BS2

| PHASE NUMBER | PERICENTER (FGS) | PERICENTER (OPTIMISED) | RESONANCE n:m (FGS) | RESONANCE n:m (OPTIMISED) |
|---|---|---|---|---|
| 1 | 6.40e5 km | 6.36e5 km | 3:1 | 3:1 |
| 2 | 6.25e5 km | 6.25e5 km | 5:2 | 5:2 |
| 3 | 6.02e5 km | 6.05e5 km | 2:1 | 2:1 |
| 4 | 5.72e5 km | 5.77e5 km | 8:5 | 8:5 |

### 4.4 - OPTIMAL ASSEMBLED SOLUTION

The two parts of the trajectory, the transfer phase from the Earth to Jupiter and the tour phase, were assembled together to form a single NLP problem for a single trajectory. The dynamical model takes into account the gravitational perturbations due to the Sun when the spacecraft is at Jupiter, and the gravity perturbation due to Jupiter during the approach phase to the giant planet. The resulting optimal solution is represented in Figure 13 for the transfer phase and in Figure 14 for the Jovian tour. Thrust arcs are represented with solid lines and coast arcs with dashed lines.

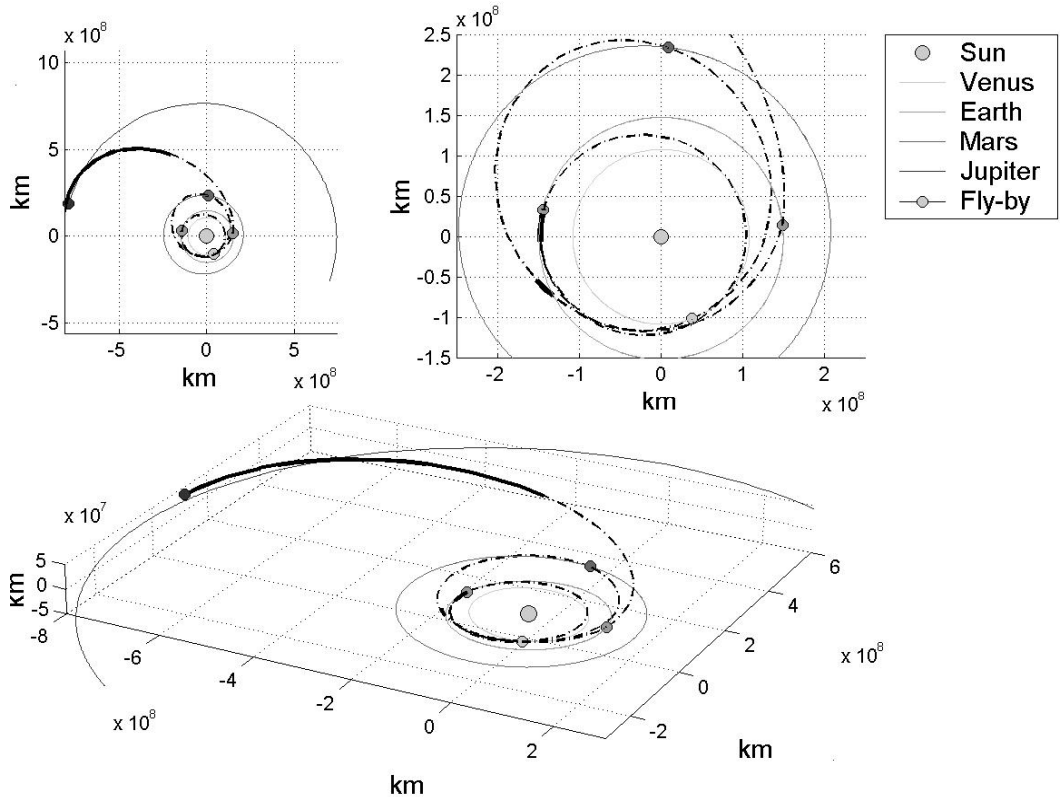

Figure 13. Optimised solution : interplanetary transfer.

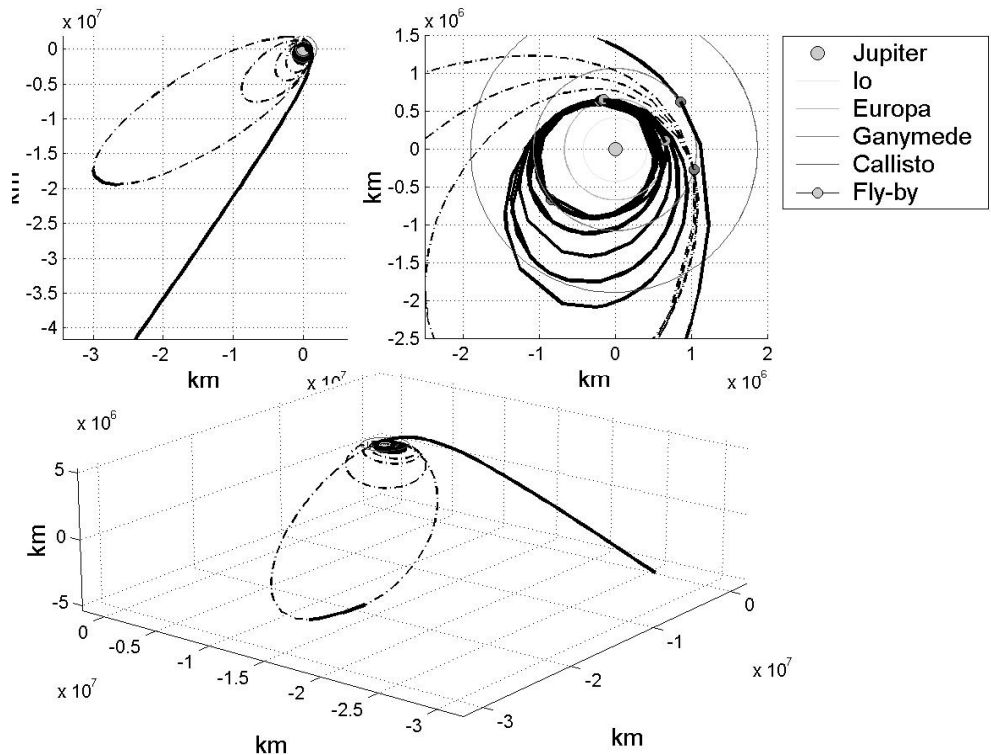

Figure 14. Optimised solution : Jovian system.

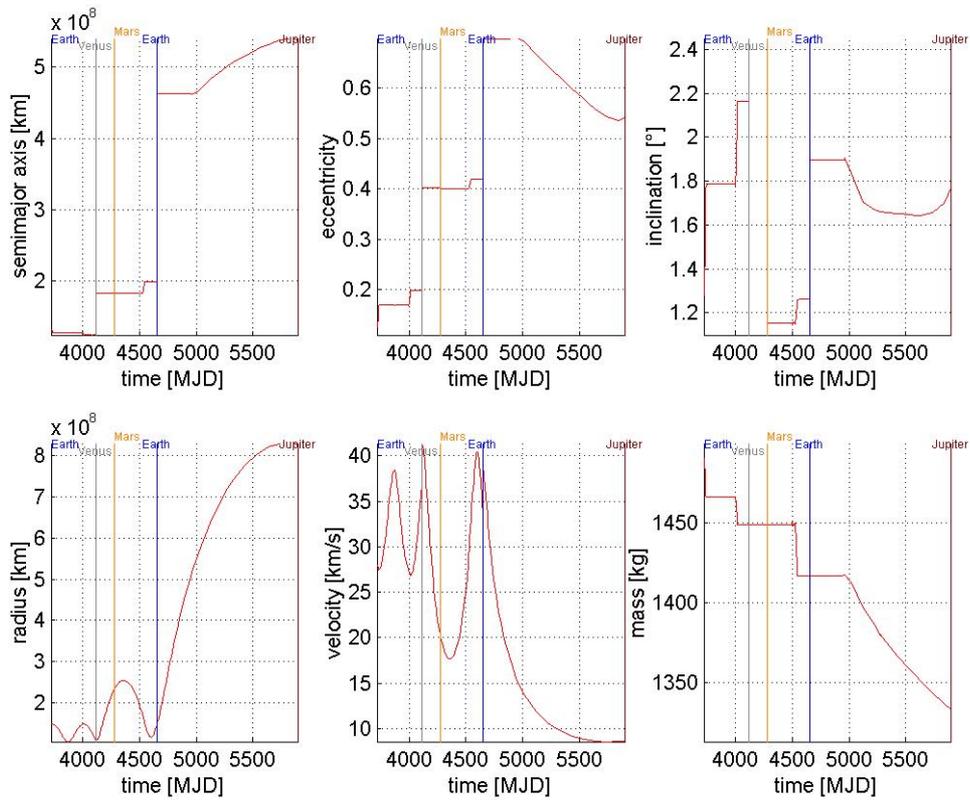

Figure 15. Orbital parameters during the interplanetary transfer

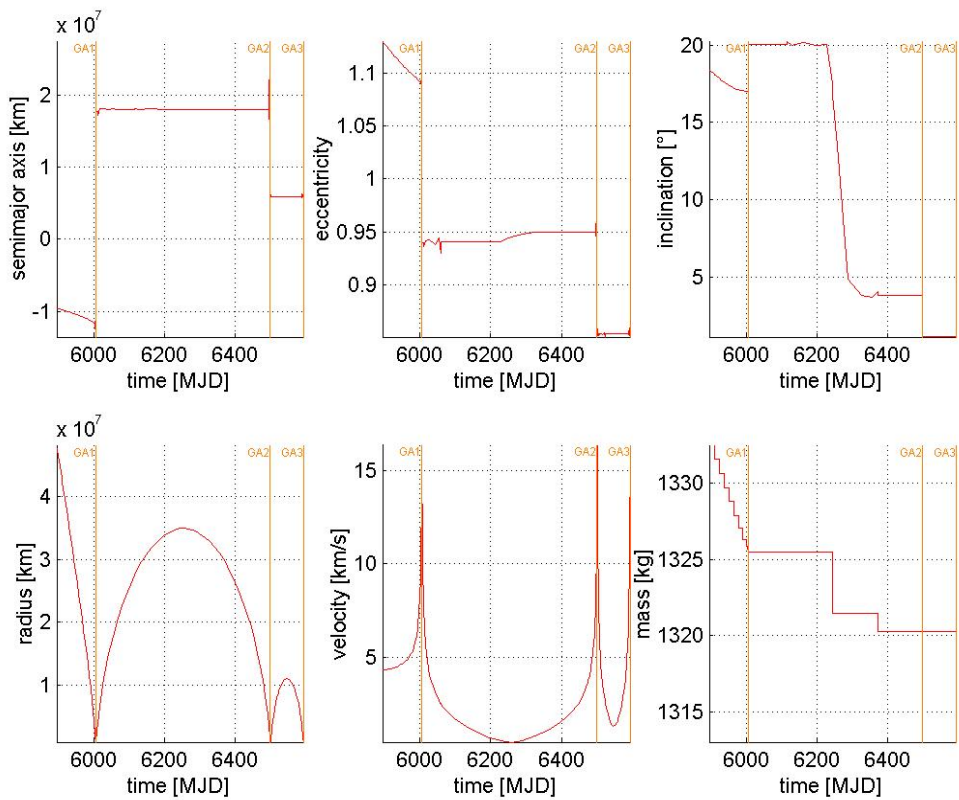

Figure 16. Orbital parameters during capture (from Ganymede flyby 1 to 3).

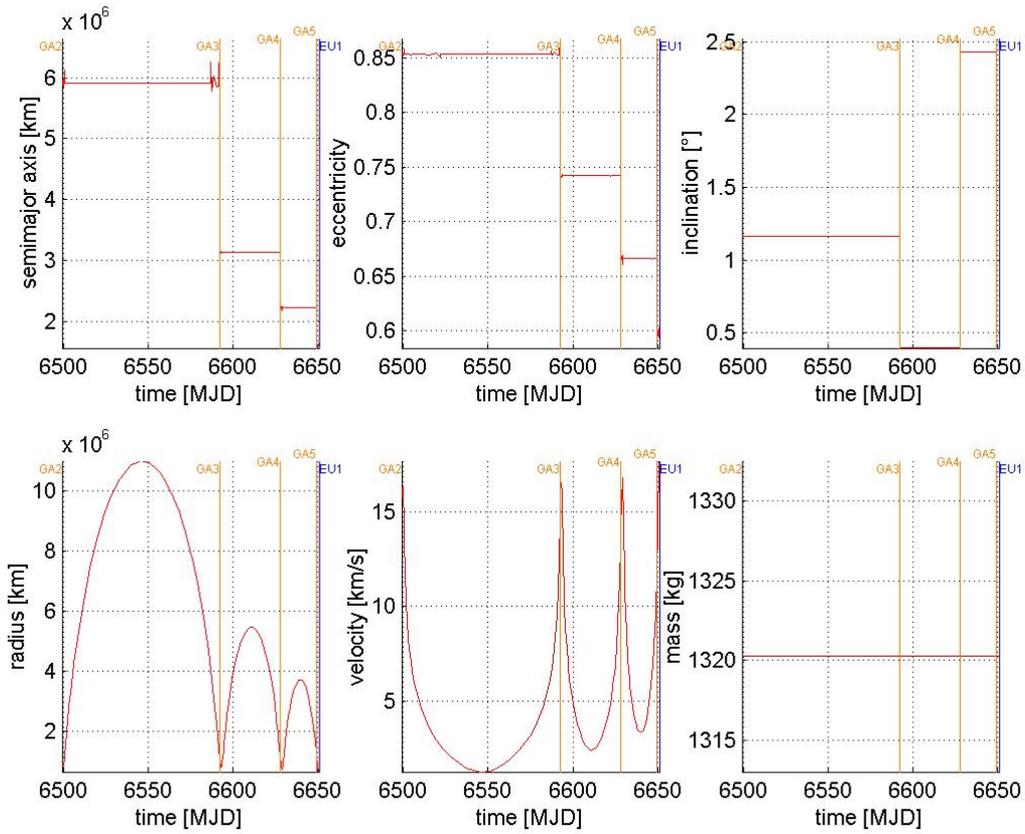

Figure 17. Orbital elements with respect to Jupiter, distance and velocity with respect to Ganymede and mass during the synchronous tour of Ganymede.

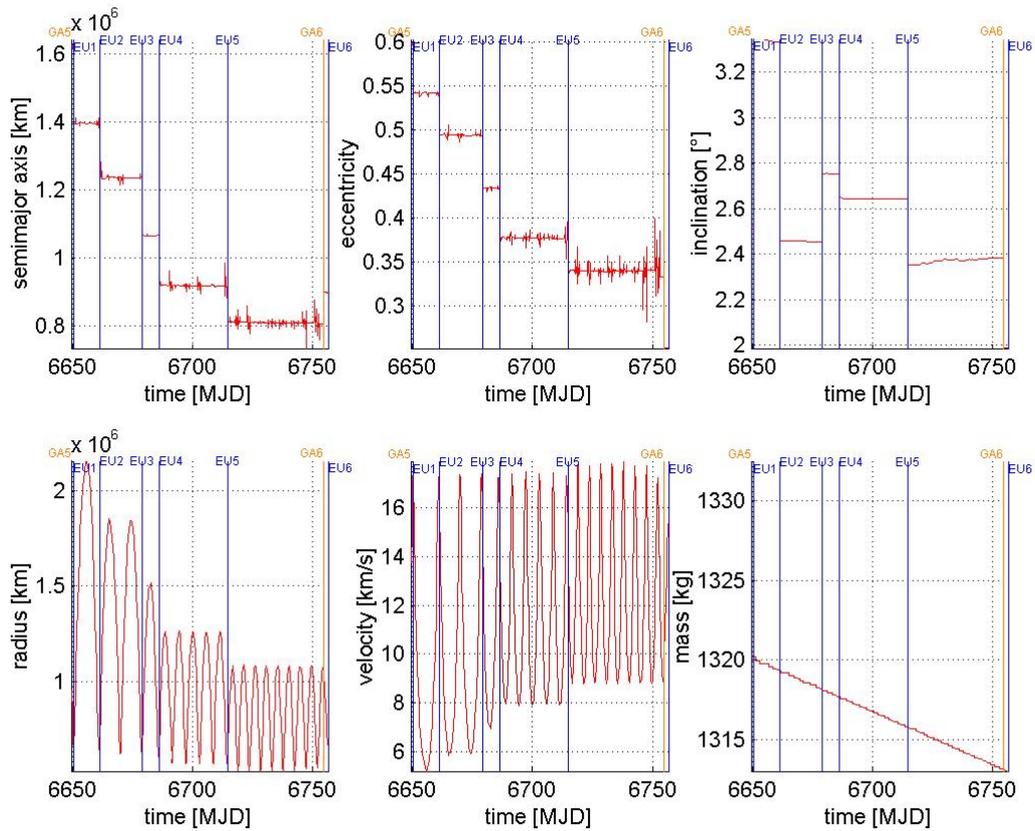

Figure 18. Orbital elements with respect to Jupiter, distance and velocity with respect to Ganymede and mass during the synchronous tour of Europae.

Figure 15 to Figure 18, represent the time history of semimajor-axis, inclination, eccentricity, velocity modus, radius and mass for the transfer trajectory, the capture and the complete tour, showing the effect of each swing-by while Table 7 gives the complete time history and the times of flight (TOF) for the entire optimised transfer. The interplanetary part of the transfer corresponds to the solution in Table 3, however when assembling all the phases together, adding also the Jupiter-centred part of the transfer, the optimiser changes a bit the dates of departure and of encounter with the planets and the propellant mass consumption for the interplanetary leg increases to match the required velocity at Jupiter to have a capture. Table 7 is a summary of the characteristics of the gravity assist manoeuvres for the entire optimised transfer.

Table 7. Time of Flights and encouter dates for the entire optimised solution

| PH N. | CELESTIAL BODY | DEP. DATE (MJD) | TOF (D) | CUMMULATIVE MISSION TIME(Y) |
|---|---|---|---|---|
| 1 | Earth | 3718.9 | 393.7 | 1.08 |
| 2 | Venus | 4112.7 | 166.5 | 1.53 |
| 3 | Mars | 4279.1 | 375.3 | 2.56 |
| 4 | Earth | 4654.4 | 1241.6 | 5.96 |
| 5 | SOI of Jupiter | 5896.1 | 110.6 | 6.27 |
| 6 | Ganymede1 | 6006.7 | 492.7 | 7.62 |
| 7 | Ganymede2 | 6499.3 | 93.0 | 7.87 |
| 8 | Ganymede3 | 6592.3 | 35.8 | 7.97 |

| 9 | Ganymede4 | 6628.1 | 21.5 | 8.03 |
| 10 | Ganymede5 | 6649.6 | 1.2 | 8.03 |
| 11 | Europa1 | 6650.8 | 10.7 | 8.06 |
| 12 | Europa2 | 6661.4 | 17.8 | 8.11 |
| 13 | Europa3 | 6679.2 | 7.1 | 8.13 |
| 14 | Europa4 | 6686.3 | 28.4 | 8.21 |
| 15 | Europa5 | 6714.7 | 39.6 | 8.32 |
| 16 | Ganymede6 | 6754.4 | 2.0 | 8.32 |

Table 8. Summary of gravity assist characteristics for the entire optimised solution

| PH N. | FROM | TO | DEP. RELATIVE VELOCITY (KM/S) | DEP. ABSOLUTE VELOCITY (KM/S) | ARRIVAL RELATIVE VELOCITY (KM/S) | ARRIVAL ABSOLUTE VELOCITY (KM/S) | B [°] | PERICENTRE ALTITUDE (KM) | $\Delta V_{GA}$ (KM/S) |
|---|---|---|---|---|---|---|---|---|---|
| 1 | Earth | Venus | 1.78 | 28.32 | 6.71 | 36.89 | 61.66 | 753 | 6.88 |
| 2 | Venus | Mars | 6.71 | 41.31 | 5.63 | 20.08 | 6.14 | 17478 | 0.60 |
| 3 | Mars | Earth | 5.63 | 20.10 | 12.26 | 33.19 | 33.09 | 300 | 6.98 |
| 4 | Earth | SOI of Jupiter | 12.26 | 38.52 | 4.30 | 8.63 | - | - | - |
| 5 | SOI of Jupiter | Ga1 | - | 4.30 | 6.40 | 15.74 | 9.01 | 200 | 1.01 |
| 6 | Ga1 | Ga2 | 6.40 | 15.16 | 6.60 | 15.17 | 8.46 | 217 | 0.97 |
| 7 | Ga2 | Ga3 | 6.60 | 14.69 | 6.61 | 14.68 | 8.39 | 235 | 0.97 |
| 8 | Ga3 | Ga4 | 6.61 | 14.02 | 6.61 | 14.02 | 8.49 | 200 | 0.98 |
| 9 | Ga4 | Ga5 | 6.61 | 13.42 | 6.60 | 13.42 | 8.50 | 200 | 0.98 |
| 10 | Ga5 | Eu1 | 6.60 | 12.60 | 4.25 | 17.26 | 10.48 | 200 | 0.78 |
| 11 | Eu1 | Eu2 | 4.25 | 16.89 | 4.29 | 16.88 | 8.89 | 512 | 0.66 |
| 12 | Eu2 | Eu3 | 4.29 | 16.53 | 4.35 | 16.50 | 9.25 | 366 | 0.70 |
| 13 | Eu3 | Eu4 | 4.35 | 16.00 | 4.37 | 15.99 | 9.61 | 268 | 0.73 |
| 14 | Eu4 | Eu5 | 4.37 | 15.39 | 4.45 | 15.36 | 8.83 | 379 | 0.68 |
| 15 | Eu5 | Ga6 | 4.45 | 14.76 | 2.04 | 8.90 | 39.10 | 2105 | 1.36 |
| 16 | Ga6 | Eu6 | 2.04 | 9.77 | 1.65 | 15.36 | - | - | - |

## 5 - CONCLUSIONS

In this paper, a mission to Europa was designed considering solar electric propulsion as main source of thrust. We devised a number of optimisation tools to design the whole trajectory: three deterministic algorithms to find a global optimum for multiple gravity assist trajectories and a general optimisation tool for optimal trajectory design based on finite elements in time. The global search was extremely efficient at providing a wide range of good first guesses. In particular, the large solution space of all the possible SOTs of Jovian moons was explored in a few seconds. Despite the fact that the first guesses were only optimal for a multi-impulse trajectory model, they were sufficiently good for the initialisation of the optimisation with DFETs. The resulting solution shows how a mission to Europa using solar electric propulsion could be feasible and interesting. On the other hand, several problems are still open and deserve further investigation, in particular, the total dose of radiation for the designed tour. But many other efficient tours can be easily designed. Also the final capture into a stable orbit around Europa, which still requires a chemical manoeuvre, could be performed with low-thrust propulsion.

## 6 - AKNOWLEDGEMENTS

The work in this paper was performed in 2002 when Massimiliano Vasile was part of the ESA, Advanced Concept Team at ESTEC and Stefano Campagnola was a student at Politecnico di Milano and an intern at the Advance Concept Team. The paper was presented in the same year at the 2nd Low Thrust International Symposium, CNES, IAS, Toulouse, France, 18-20 June 2002.

This work was funded in part by the European Space Agency. The authors would like to thank Dr. Guy Janin, Mr. Franco Ongaro, Dr. Saccoccia and Dr. Robin Biesbroek at the European Space Agency for their generous support and expert guidance.

## 7 - BIBLIOGRAPHY